\documentclass [a4paper,12pt]{amsart}
\usepackage{graphicx}

\usepackage{amssymb}

\newtheorem{thm}{Theorem}[section]
\newtheorem{cor}[thm]{Corollary}
\newtheorem{lem}[thm]{Lemma}

\theoremstyle{definition}

\theoremstyle{definition}
\newtheorem{rem}[thm]{Remark}

\theoremstyle{definition}

\def\R{\mathbb R}
\def\re{\mathbb R}

\def\sign{\operatorname{sign}}

\def\ind{\operatorname{Ind}}

\begin{document}

\author{X. Gual-Arnau and J.J.~Nu\~no-Ballesteros}

\title{A stereological version of the Gauss-Bonnet formula}


\date{}

\address{X. Gual-Arnau, Department de Matem\`atiques, Campus
Riu Sec, Universitat  Jaume I, 12071 Castell\'o SPAIN}

\email{gual@mat.uji.es}

\address{J.J.~Nu\~no-Ballesteros, Departament de Geometria i Topologia,
Universitat de Val\`encia, Campus de Burjassot, 46100 Burjassot
SPAIN}

\email{nuno@uv.es}

\thanks{The first author is partially supported by  Generalitat Valenciana Grant
GV-97-CB-12-72. The second author is partially supported by DGICYT
Grant PB96--0785}

\maketitle

\section{Introduction}
Let $D\subset S\subset \re^3$ be a domain with boundary in an
orientable smooth surface $S$ in $\re^3$. Then the Gauss-Bonnet
formula gives a method to compute the Euler-Poincar\'e
characteristic of $D$ in terms of the Gauss curvature $K$ of $S$,
and the geodesic curvature $k_g$ of $\partial D$ in $S$:
$$2\pi\chi(D)=\int_D K+\int_{\partial D} k_g.$$ This can be found
in any elementary book of Classical Differential Geometry. In this
paper, we give a stereological version of this formula, in the
following sense: let $u\in S^2$ and denote by $\pi_{u,\lambda}$
the plane orthogonal to $u$ given by the equation $\langle x,
u\rangle=\lambda$.  When $\lambda $ varies on $\R$, the different
planes $\pi_{u,\lambda}$ can be considered as a ``sweeping'' plane
in $\R^3$ and its contact with $D$ and with $\partial D$ will give
a method to compute $\chi(D)$: $$\chi(D)=\sum_{x\in A} \sign
K(x)+\frac 12\sum_{y\in B} \sign (k_g(y)-k_g^u(y)).$$ Here, $A$ is
given by the points $x\in D$ such that $D$ is tangent to
$\pi_{u,\lambda}$ at $x$ for some $\lambda$, $B$ is equal to the
subset of points $y\in\partial D$ where $\partial D$ is tangent to
$\pi_{u,\lambda}\cap S$ for some $\lambda$, and $k_g^u$ denotes
the geodesic curvature of $\pi_{u,\lambda}\cap S$. We will show
that for a generic $u\in S^2$:
\begin{enumerate}
\item the sets $A,B$ are finite;
\item if $x\in A$, then $K(x)\ne 0$;
\item if $y\in B$, then $\pi_{u,\lambda}\cap S$ is a regular curve
in $S$ and $k_g(y)\ne k_g^u(y)$.
\end{enumerate}
In this way, the above formula has sense for almost any $u\in S^2$
(that is, for any $u\in S^2\smallsetminus N$, where $N$ is a null
set in $S^2$). Moreover, this generalizes the result of
\cite{GBNO}, which computes the Euler-Poincar\'e characteristic of
a plane domain $D\subset \R^2$ by looking at contact with lines.

\section{Proof of the formula}

The classical Poincar\'e-Hopf Theorem stays that if $S$ is a
closed orientable smooth surface and $v$ is a smooth vector field
on $S$ with isolated zeros, then
$$\chi(S)=\sum_{v(x)=0}\ind_x(v),$$ where $\chi(S)$ is the
Euler-Poincar\'e characteristic of $S$ and the index, $\ind_x(v)$,
is just the local degree of $v$ at $x$. As an immediate
consequence, we get that if $f:S\to\re$ is a Morse function (that
is, a function with non-degenerate critical points), then
$$\chi(S)=\sum_{x\in\Sigma(f)}\ind_x(f),$$ where $\Sigma(f)$
denotes the set of singular points of $f$ and the index,
$\ind_x(f)$, is given by $+1$ when $x$ is a local extreme of $f$,
or $-1$ when $x$ is a saddle point.

The Poincar\'e-Hopf Theorem was generalized for surfaces with
boundary by Morse \cite{Mo} in the following way: suppose that $D$
is a compact orientable smooth surface with boundary (that we can
assume embedded in an orientable smooth surface $S$ without
boundary) and let $v$ be a smooth vector field on $D$ such that
\begin{enumerate}
\item $v$ has isolated zeros on $D$,
\item $v$ is not zero on $\partial D$, and
\item $v$ is tangent to $\partial D$ only at a finite number of
points. \end{enumerate}
If $y\in \partial D$ is a point where $v$
is tangent to $\partial D$, generically we have that the integral
curve of $v$ at $y$ is locally contained in $S\smallsetminus D$ or
in $D$ and we can assign an index $+1$ or $-1$ respectively. In
the general case, we can extend this by using the Transversality
Theorem and define an index, $\ind_y(v)\in\{-1,0,+1\}$. Then, the
Morse Theorem stays that $$\chi(D)=\sum_{v(x)=0}\ind_x(v)+\frac
12\sum_{v(y)\in T_y\partial D}\ind_y(v).$$

Now, we can interpret this in terms of critical points of
functions. Suppose that $f:D\to \re$ is a Morse function such that
$f$ has no critical points on $\partial D$ and the restriction
$f|_{\partial D}:\partial D\to \re$ is also a Morse function. For
each point $y\in\Sigma(f|_{\partial D})$, we have that the level
set $f^{-1}(f(y))$ is locally contained either in $S\smallsetminus
D$ or in $D$. In the first case, we say that $y$ is a point of
type ``island'' and we assign an index $\ind_{y}(f)=+1$ and in the
second case, we say that it is point of type ``bridge'' and assign
an index $\ind_{y}(f)=-1$ (see Figure 1).

\bigskip
\centerline{\includegraphics{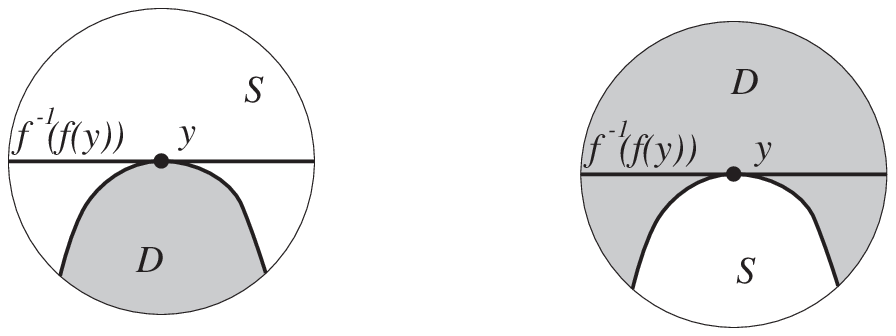}}

\medskip
\centerline{a) ``Island'' \hskip45mm b) ``Bridge''}
\medskip
\centerline{Figure 1}
\bigskip

\begin{lem} Let $D$ be a compact orientable smooth surface
with boundary and let $f:D\to \re$ be a Morse function such that
$f$ has no critical points on $\partial D$ and the restriction
$f|_{\partial D}:\partial D\to \re$ is also a Morse function.
Then, $$\chi(D)=\sum_{x\in\Sigma(f)} \ind_{x}(f)+\frac 12
\sum_{y\in\Sigma(f|_{\partial D})} \ind_{y}(f).$$
\end{lem}

An interesting application of this formula is obtained when
$D\subset S\subset \re^3$ and we consider the family of height
functions. Given $u\in S^2$, we define the height function
$h_u:\re^3\to \re$ by $h_u(x)=\langle x,u\rangle$. When we
restrict this function to $D$ (or $S$), the level sets are the
plane curves $\pi_{u,\lambda}\cap D$ (resp. $\pi_{u,\lambda}\cap
S$), where $\pi_{u,\lambda}$ is the plane $h_u^{-1}(\lambda)$. In
particular, if $\lambda$ is a regular value of $h_u|_D$ (and hence
of $h_u|_S$), we have that $\pi_{u,\lambda}\cap S$ is in fact a
smooth curve in $S$.

We will give a geometrical interpretation of the critical points
of $h_u|_D$ and $h_u|_{\partial D}$ and the corresponding indices.
Before that, we need to fix some notation. Given $x\in D$, we
denote by $N(x)$ and $K(x)$ the normal vector and the Gauss
curvature of $D$ at $x$ respectively. Given $y\in
\partial D$, we denote by $n(y)$ and $k_g(y)$ the normal vector
and the geodesic curvature of $\partial D$ in $D$ at $y$
respectively (we choose the orientation in $\partial D$ so that
$n(y)$ points to the interior of $D$). Finally, we denote by
$k_g^u(y)$ the geodesic curvature of the curve
$\pi_{u,\lambda}\cap S$ at a regular point $y$ of $h_u|_D$ (in the
case that $\pi_{u,\lambda}\cap S$ and $\partial D$ are tangent at
$y$, we choose in $\pi_{u,\lambda}\cap S$ the same orientation).

\begin{thm} Let $D\subset S\subset \re^3$ be a domain with boundary in
an orientable smooth surface $S$ in $\re^3$ and let $u\in S^2$.
\begin{enumerate}
\item Let $x\in D$. It is a critical point of $h_u|_D$ if and only if
$\pi_{u,\lambda}$ is tangent to $D$ at $x$ (where
$h_u(x)=\lambda$) if and only if $u=\pm N(x)$.

\item Let $x\in D$ be a critical point of $h_u|_D$. It
is non-degenerate if and only if $K(x)\ne 0$. Moreover, it is a
local extreme of $h_u|_D$ when $K(x)>0$ and a saddle point when
$K(x)>0$.

\item Let $y\in \partial D$ be a regular point of $h_u|_D$.
It is a critical point of $h_u|_{\partial D}$ if and only if
$\pi_{u,\lambda}\cap S$ is tangent to $\partial D$ at $y$ if and
only if $\tilde u=\pm n(y)$, where $\tilde u$ denotes the
normalized orthogonal projection of $u$ in $T_yD$.

\item Let $y\in \partial D$ be a regular point of $h_u|_D$ which
is also a critical point of $h_u|_{\partial D}$. It is
non-degenerate if and only if $k_g(y)\ne k_g^u(y)$. Moreover, it
is an island when $k_g(y)>k_g^u(y)$ and a bridge when
$k_g(y)<k_g^u(y)$.
\end{enumerate}
\end{thm}

\begin{proof} Parts (1) and (2) are standard results on Generic Geometry
and can be found for instance in \cite{BG}. We will prove (3) and
(4). Suppose that $\alpha(s)$ is a parametrization of $\partial D$
by arc length, so that $y=\alpha(s_0)$ and
$n(y)=N(y)\wedge\alpha'(s_0)$. Since $\{N(y),\alpha'(s_0),n(y)\}$
is a orthonormal basis of $\re^3$, we have that $$u=\langle
u,N(y)\rangle N(y)+\langle u,\alpha'(s_0)\rangle
\alpha'(s_0)+\langle u,n(y)\rangle n(y).$$ In particular, provided
that $u\ne\pm N(y)$, we have that $\tilde u=v/\|v\|$, where
$$v=\langle u,\alpha'(s_0)\rangle \alpha'(s_0)+\langle
u,n(y)\rangle n(y).$$ Thus, $\tilde u=\pm n(y)$ if and only if
$\langle u,\alpha'(s_0)\rangle=0$ if and only if $y$ is a critical
point of $h_u|_{\partial D}$.

On the other hand, let $\beta(s)$ be a parametrization of
$\pi_{u,\lambda}\cap S$ by arc length, with $\beta(s_0)=y$. In
this case, $\langle u,\beta(s)\rangle=\lambda$ for all $s$ and
thus $\langle u,\beta'(s_0)\rangle=0$. This gives that $\tilde u$
is also equal to the normal vector to $\beta'(s_0)$ in
$T_y\partial D$. In particular, $y$ is a critical point of
$h_u|_{\partial D}$ if and only if the two curves are tangent.

\medskip
Suppose now that $y$ is a critical point of $h_u|_{\partial D}$.
We have that $$\alpha''(s_0)=k_n(y)N(y)+k_g(y)n(y),$$ where
$k_n(y)$ is the normal curvature of $\alpha$. This gives that
$$\langle u,\alpha''(s_0)\rangle=k_n(y)\langle u,N(y)\rangle
+k_g(y)\langle u,n(y)\rangle.$$

But we can do the same with $\beta(s)$, getting $$0=\langle
u,\beta''(s_0)\rangle=k_n(y)\langle u,N(y)\rangle +k_g^u(y)\langle
u,n(y)\rangle.$$ This implies that $$\langle
u,\alpha''(s_0)\rangle=(k_g(y)-k_g^u(y))\langle u,n(y)\rangle.$$
Since $\langle u,n(y)\rangle\ne 0$, we conclude that $y$ is a
non-degenerate critical point of $h_u|_{\partial D}$ if and only
if $k_g(y)\ne k_g^u(y)$. The fact that the cases $k_g(y)>
k_g^u(y)$ and $k_g(y)< k_g^u(y)$ correspond to an island or a
bridge respectively can be deduced from the way we have chosen the
orientation on $\partial D$.
\end{proof}

The following consequence can be seen as a stereological version
of the Gauss-Bonnet formula.

\begin{cor}\label{espai} Let $D\subset S\subset \re^3$ be a domain with boundary in
an orientable smooth surface $S$ in $\re^3$. For a generic $u\in
S^2$ the height function $h_u|_D$ is a Morse function, has no
critical points on $\partial D$ and its restriction to $\partial
D$ is also a Morse function. In particular, we have that
$$\chi(D)=\sum_{x\in D/N(x)=\pm u} \sign K(x)+\frac
12\sum_{y\in\partial D/n(y)=\pm \tilde u} \sign
(k_g(y)-k_g^u(y)).$$
\end{cor}

\begin{proof} The formula for $\chi(D)$ is a direct consequence of Lemma
3.1 and Theorem 3.2. We will see the first part of the statement.
Given $u\in S^2$, $h_u|_D$ is a Morse function if and only if $u$
is a regular value of the Gauss map $N:D\to S^2$ (just note that
the points $x\in D$ where $K(x)=0$ are the critical points of the
Gauss map). Therefore, the Sard Theorem implies that for a generic
$u\in S^2$, $h_u|_D$ is a Morse function. The other two conditions
are also an easy consequence of this theorem.
\end{proof}

It is also interesting to look at some particular cases. Suppose
for instance, that $S$ is a plane, namely $S=\re^2$. Then we can
use height functions $h_u:\re^2\to \re$ for $u\in S^1$, instead of
the height functions of $\re^3$. The level curves in this case are
the orthogonal lines to $u$, which have zero curvature. Moreover,
$h_u$ has no critical points in $\re^2$ and the geodesic curvature
is equal to the classical curvature of a plane curve. Thus, we get
the following formula method to obtain the Euler number of a
domain in $\re^2$ based on 1-dimensional observations, modified
from \cite{H}, and that has been used in Stereology, \cite{GBNO}.

\begin{cor}\label{pla} Let $D\subset \re^2$ be a domain with boundary. For a
generic $u\in S^1$ the restriction of the height function
$h_u|_{\partial D}$ is a Morse function and $$\chi(D)=\frac 12
\sum_{y\in\partial D/n(y)=\pm u} \sign k(y).$$
\end{cor}

In the case that $S=S^2$, the level curves associated to $u\in
S^2$ are the corresponding parallels. Note that $h_u|_{S^2}$ is a
Morse function with two critical points, $u$ and $-u$. Moreover,
if we denote by $\theta_u,\gamma_u$ the spherical coordinates
associated to $u$, we have that the geodesic curvature of the
parallel $\gamma_u=\gamma_0$ is $-\tan\gamma_0$.

\begin{cor} Let $D\subset S^2$ be a domain with boundary. For a generic
$u\in S^2$, $$\chi(D)=\frac 12 \sum_{y\in\partial D/n(y)=\pm
\tilde u} \sign (k_g(y)+\tan\gamma_u(y))+\#(\{u,-u\}\cap D),$$
where $\#(\{u,-u\}\cap D)$ means the number of times that $u$ or
$-u$ belong to $D$.
\end{cor}

It is also possible to compute the Euler number of $D\subset S^2$
from observations on a meridian that sweeps through $D$. This is
not a consequence of the construction associated to the height
functions, but a different family of functions. Let $u\in S^2$ and
let $S^1$ be the circle in $S^2$ orthogonal to $u$. Then, we have
the regular function $\theta_u:S^2\smallsetminus\{u,-u\}\to S^1$,
whose level sets are the meridians $\theta_u=\theta_0$ (which have
zero geodesic curvature). Moreover, for a generic $u\in S^2$, the
restriction of $\theta_u$ to $\partial D$ is a Morse function,
\cite{T}.

\begin{cor} Let $D\subset S^2$ be a domain with boundary. For a generic
$u\in S^2$, $$\chi(D)=\frac 12 \sum_{y\in\partial D/n(y)\bot
\tilde u} \sign k_g(y)+\#(\{u,-u\}\cap D).$$
\end{cor}

\begin{proof} Let $S=S^2\smallsetminus\{u,-u\}$ and let $\tilde D$ be the
domain with boundary in $S$ given by $D\smallsetminus (B_1\cup
B_2)$, where $B_1,B_2$ are small geodesic balls centered at $u,-u$
which do not intersect $\partial D$. It follows from Lemma 3.1
that $$\chi(\tilde D)=\frac 12 \sum_{y\in\partial D/n(y)\bot
\tilde u} \sign k_g(y).$$ On the other hand, $\chi(D)=\chi(\tilde
D)+\#(\{u,-u\}\cap D)$, which concludes the proof. \end{proof}

\begin{rem} The preceding result may be extended to ovaloids in the
following way: Let $S$ be an ovaloid (compact and connected
surface in $\re^3$ for which the Gaussian curvature $K>0$) and
$\tilde{D}$ a domain in $S$; then, $S$ is diffeomorphic to a
sphere $S^2$ through its Gauss map. Let $D\subset S^2$ denote the
image of $\tilde{D}$ under this diffeomorphism; then, $\chi
(D)=\chi (\tilde{D})$ and it is possible to obtain $\chi
(\tilde{D})$ from observations in the curves which are image of
the meridians in $S^2$. \end{rem}

\section{Stereological applications}

The Euler-Poincar\'e characteristic for domains in $\R^2$ and
$\R^3$ has been studied in several stereological applications
(see, for instance \cite{GBNO}). The principle used to obtain the
Euler number of an $n$-dimensional domain in $\R^n$ ($n=2$ or
$3$), is based on what happens in an $(n-1)$-dimensional plane
that sweeps through the domain. For domains in $\R^2$, this
principle is based on Corollary \ref{pla}. Now, we will extend
this principle to domains in a surface.

Let $D\subset S\subset \R^3$ be a domain with boundary in an
orientable smooth surface $S$ in $\R^3$ and let $u\in S^2$. When
$\lambda $ varies on $\R$, the different planes $\pi_{u,\lambda}$
can be considered as a ``sweeping'' plane in $\R^3$, and Corollary
\ref{espai} can be expressed as: $$\chi(D)=(I_2-B_2)+\frac
12(I_1-B_1),$$ where $I_1,B_1$ denote the number of islands and
bridges, respectively, observed in the level curves
$\pi_{u,\lambda}\cap S$ (see Figure 2) and $I_2,B_2$, similarly to
\cite{GBNO}, denote the number of islands and bridges observed in
the ``sweeping'' plane $\pi_{u,\lambda}$, which contribute to the
sum $\sum_{x\in\Sigma(h_u)}\ind_x(h_u)$ (see Figure 2).

\bigskip
\centerline{\includegraphics{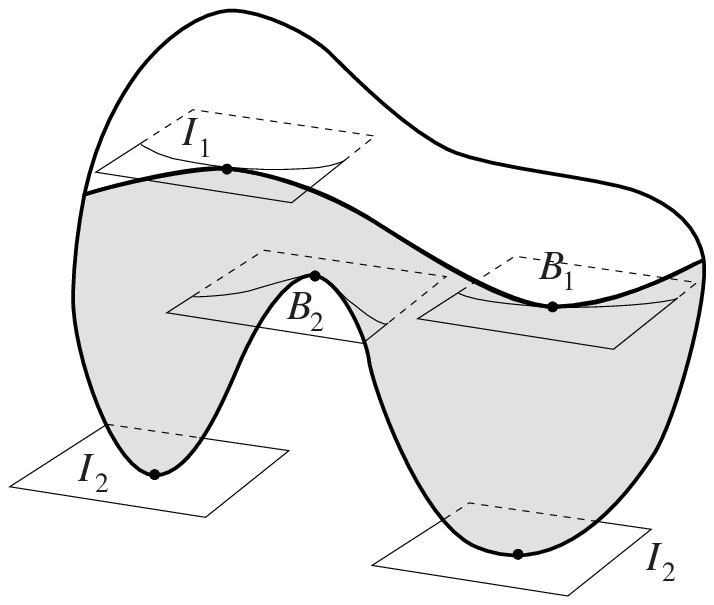}}

\medskip
\centerline{Figure 2}


\bigskip


\begin{thebibliography}{1}


\bibitem{BG}
J.W. Bruce, P.J. Giblin, \emph{Curves and singularities,}
Cambridge University Press, 1984.


\bibitem{GBNO}
H. J. G. Gundersen, R. W. Boyce, J.R. Nyengaard, A. Odgaard,
\emph{The conneulor: unbiased estimation of connectivity using
physical disectors under projection,} Bone \textbf{14} (1993),
217--222.

\bibitem{H}
H. Hadwiger, \emph{Vorlesungen, \"Uber Inhalt, Oberfl\"ache und
Isoperimetrie,} Springer-Verlag, Berlin, 1957.


\bibitem{Mo}
M. Morse, \emph{Singular points of vector fields under general
boundary conditions,} American Journal of Mathematics \textbf{51}
(1929), 165--178.


\bibitem{S}
L. A. Santal\'o, \emph{Integral Geometry and Geometric
Probability,} Addison-Wesley Publishing Company Inc., London,
1976.

\bibitem{T} E. Teufel, \emph{Differential Topology and the
Computation of Total Absolute Curvature,}  Math. Ann  \textbf{258}
(1982), 471--480


\end{thebibliography}
\end{document}